\RequirePackage{fix-cm} % Fix LaTeX2e bugs.
\documentclass[a4paper, twoside, reqno, dvips, 12pt]{amsart}
\usepackage{fixltx2e}   % Fix LaTeX2e bugs.

\usepackage[latin1]{inputenc}
\usepackage[T1]{fontenc}

\usepackage{eucal}
\usepackage{esint}
\usepackage{dsfont}
\usepackage{xspace}
\usepackage{amsgen}
\usepackage{amsthm}
\usepackage{amssymb}
\usepackage{amsmath}
\usepackage{upgreek}
\usepackage{amsfonts}
\usepackage{MnSymbol}
\usepackage{mathrsfs}
\usepackage{mathtools}
\usepackage[nice]{nicefrac}

\usepackage{a4wide}

\usepackage[dvipsnames, table]{xcolor}

\definecolor{oneblue}{rgb}{0.0, 0.0, 0.85}
\definecolor{bluepigment}{rgb}{0.2, 0.2, 0.6}
\definecolor{darkgrey}{rgb}{0.273, 0.281, 0.30}
\definecolor{Lightgray}{rgb}{0.89, 0.89, 0.89}
\definecolor{Lightblue}{RGB}{214, 214, 214}
\definecolor{bckg}{RGB}{20.8, 20.8, 20.8} % Color of the boxes
\definecolor{charcoal}{rgb}{0.21, 0.27, 0.31}
\definecolor{darkelectricblue}{rgb}{0.33, 0.41, 0.47}

\usepackage[sort&compress, comma, square, numbers]{natbib}

\usepackage{psfrag}
\usepackage{graphicx}
\usepackage{subfigure}
\usepackage{morefloats}
\usepackage{indentfirst}

\usepackage{acronym}
\usepackage{microtype}
\usepackage[labelsep=period,%
            labelfont={bf,sf,color=bluepigment},%
            justification=raggedright]{caption}

\usepackage[perpage, symbol]{footmisc}

\usepackage[usenames, dvipsnames, pdf]{pstricks}
\usepackage{epsfig}
\usepackage{pst-grad} % For gradients
\usepackage{pst-plot} % For axes

\usepackage[colorlinks,
           urlcolor=oneblue,
           linkcolor=oneblue,
           citecolor=NavyBlue,
           bookmarksopen=false,
           pdfpagemode=UseNone,
           pagebackref]{hyperref}

\usepackage[explicit]{titlesec}

\titleformat{\section}
  {\color{NavyBlue}\Large\sffamily\bfseries}
  {}
  {0em}
  {\colorbox{bckg!5}{\parbox{\dimexpr\linewidth-2\fboxsep\relax}{\centering\thesection. #1}}}
  [\vspace*{0.33em}]

\titleformat{name=\section,numberless}
  {\color{NavyBlue}\Large\sffamily\bfseries}
  {}
  {0.0em}
  {\colorbox{bckg!10}{\parbox{\dimexpr\linewidth-2\fboxsep\relax}{\centering#1}}}
  [\vspace*{0.33em}]

\titleformat{\subsection}
  {\color{NavyBlue}\large\sffamily\bfseries}
  {}
  {0.0em}
  {\colorbox{bckg!5}{\parbox{\dimexpr\linewidth-2\fboxsep\relax}{\centering\thesubsection. #1}}}
  [\vspace*{0.33em}]

\titleformat{name=\subsection,numberless}
  {\color{NavyBlue}\Large\sffamily\bfseries}
  {}
  {0em}
  {\colorbox{bckg!10}{\parbox{\dimexpr\linewidth-2\fboxsep\relax}{\centering#1}}}
  [\vspace*{0.33em}]

\titleformat{\subsubsection}
  {\color{bluepigment}\sffamily\normalsize\bfseries}
  {\thesubsubsection}
  {0.5em}
  {#1}
  [\vspace*{0.33em}]

\titleformat{\paragraph}[runin]
  {\color{bluepigment}\sffamily\small\bfseries}
  {}
  {0em}
  {#1}

\titlespacing{\section}{1.0em}{1.5em plus 2pt minus 2pt}%
{1.0em plus 2pt minus 2pt}[0em]
\titlespacing{\subsection}{1.0em}{1.5em plus 2pt minus 2pt}%
{1.0em}[0em]
\titlespacing{\subsubsection}{1.0em}{1.5em plus 2pt minus 2pt}%
{1.0em plus 2pt minus 2pt}[0em]

\usepackage{titletoc}

\contentsmargin{0.5em}
\newlength{\tocsep} 
\titlecontents{section}[\tocsep]
  {\addvspace{10pt}\bfseries\sffamily}
  {\contentslabel[\thecontentslabel]{\tocsep}}
  {}
  {\ \titlerule*[0.75pc]{.}\ \thecontentspage}
  []
\titlecontents{subsection}[\tocsep]
  {\addvspace{8pt}\sffamily}
  {\contentslabel[\thecontentslabel]{\tocsep}}
  {}
  {\ \titlerule*[0.5pc]{.}\ \thecontentspage}
  []
\titlecontents*{subsubsection}[\tocsep]
  {\addvspace{2pt}\footnotesize\sffamily}
  {}
  {}
  {\ \titlerule*[0.35pc]{.}\ \thecontentspage}
  [\\*]

\makeatletter
\def\@setauthors{%
  \begingroup
  \def\thanks{\protect\thanks@warning}%
  \trivlist
  \centering\footnotesize \@topsep30\p@\relax
  \advance\@topsep by -\baselineskip
  \item\relax
  \author@andify\authors
  \def\\{\protect\linebreak}%
  \textsc{\normalsize\textcolor{darkelectricblue}{\authors}}%
  \ifx\@empty\contribs
  \else
    ,\penalty-3 \space \@setcontribs
    \@closetoccontribs
  \fi
  \endtrivlist
  \endgroup
}
\def\@settitle{\begin{center}%
  \baselineskip14\p@\relax
    \bfseries
    \textsc{\Large\textcolor{charcoal}{\@title}}
  \end{center}%
}
\makeatother

\usepackage{enumitem}
\setlist[description]{%
  topsep=30pt,               % space before start / after end of list
  itemsep=5pt,               % space between items
  font={\bfseries\sffamily\color{NavyBlue}}, % if colour is needed
}

\usepackage{fancyhdr}
\usepackage{lastpage}

\newcommand*\Title{\textcolor{bluepigment}{Multi-Symplectic Structure for SGN equations}}
\newcommand*\Authors{\textcolor{bluepigment}{M.~Chhay, D.~Dutykh \& D.~Clamond}}
\newcommand*{\plogo}{\textcolor{gray}{{\texttt{arXiv.org} / \textsc{hal}}}} % Generic publisher logo

\numberwithin{equation}{section}

\newcommand{\vs}{\tilde{v}}

\newcommand{\udelta}{\delta}

\newcommand{\phim}{\bar{\phi}}

\newcommand{\nus}{\tilde{\nu}}

\newcommand{\ve}{\boldsymbol{e}}

\newcommand{\vz}{\boldsymbol{z}}

\newcommand{\vum}{\bar{u}}
\newcommand{\vmum}{\bar{\mu}}

\newcommand{\ie}{\emph{i.e.}}

\newcommand{\scal}{\boldsymbol{\cdot}}
\newcommand{\grad}{\boldsymbol{\nabla}}

\newcommand{\half}{{\textstyle{1\over2}}}
\newcommand{\third}{{\textstyle{1\over3}}}
\newcommand{\sixth}{{\textstyle{1\over6}}}

\begin{document}

\title[\Title]{On the multi-symplectic structure of the Serre--Green--Naghdi equations}

\author[M.~Chhay]{Marx Chhay}
\address{LOCIE, UMR 5271 CNRS, Universit\'e Savoie Mont Blanc, Campus Scientifique, 
73376 Le Bourget-du-Lac Cedex, France}
\email{Marx.Chhay@univ-savoie.fr}
\urladdr{http://marx.chhay.free.fr/}

\author[D.~Dutykh]{Denys Dutykh$^*$}
\address{LAMA, UMR 5127 CNRS, Universit\'e Savoie Mont Blanc, Campus Scientifique, 
73376 Le Bourget-du-Lac Cedex, France}
\email{Denys.Dutykh@univ-savoie.fr}
\urladdr{http://www.denys-dutykh.com/}
\thanks{$^*$ Corresponding author}

\author[D. Clamond]{Didier Clamond}
\address{Universit\'e de Nice -- Sophia Antipolis, Laboratoire J. A. Dieudonn\'e, 
Parc Valrose, 06108 Nice cedex 2, France}
\email{diderc@unice.fr}
\urladdr{http://math.unice.fr/~didierc/}

%%% ------------------------------------------------------------------------ %%%

\begin{titlepage}
\thispagestyle{empty} % Remove page numbering on this page
\noindent
{\Large Marx \textsc{Chhay}}\\
{\it\textcolor{gray}{Universit\'e Savoie Mont Blanc, Polytech Annecy--Chamb\'ery, LOCIE, France}}
\\[0.02\textheight]
{\Large Denys \textsc{Dutykh}}\\
{\it\textcolor{gray}{Universit\'e Savoie Mont Blanc, CNRS, LAMA, France}}
\\[0.02\textheight]
{\Large Didier \textsc{Clamond}}\\
{\it\textcolor{gray}{Universit\'e de Nice -- Sophia Antipolis, LJAD, France}}
\\[0.16\textheight]

\colorbox{Lightblue}{
  \parbox[t]{1.0\textwidth}{
    \centering\huge\sc
    \vspace*{0.99cm}
    
    \textcolor{bluepigment}{On the multi-symplectic structure of the Serre--Green--Naghdi equations}
    
    \vspace*{0.7cm}
  }
}

\vfill % Whitespace between the title block and the publisher

\raggedleft     % Right-align all text
{\large \plogo} % Publisher and logo
\end{titlepage}

%%% ------------------------------------------------------------------------ %%%

\newpage
\maketitle
\thispagestyle{empty}

%%% ------------------------------------------------------------------------ %%%

\begin{abstract}

In this short note, we present a multi-symplectic structure of the Serre--Green--Naghdi (SGN) equations modelling nonlinear long surface waves in shallow water. This multi-symplectic structure allow the use of efficient finite difference or pseudo-spectral numerical schemes preserving \emph{exactly} the multi-symplectic form at the discrete level.

\bigskip
\noindent \textbf{\keywordsname:} fully nonlinear long waves; Serre equations; Green--Naghdi model; multi-symplectic structure \\


\end{abstract}

%%% ------------------------------------------------------------------------ %%%

\newpage
\tableofcontents
\thispagestyle{empty}

%%% ------------------------------------------------------------------------ %%%

\newpage
\section{Introduction}

The present manuscript is devoted to a further study of the celebrated Serre--Green--Naghdi model of fully nonlinear long water waves propagating in shallow water. Namely, we unveil another variational structure of these equations. The Hamiltonian formulation for the Serre equations can be found, for example, in \cite{Johnson2002}. However, this structure is non-canonical and highly non-trivial, at least at the first sight. In this article, we propose a multi-symplectic structure for the same system of Serre equations. The multi-symplectic structure generalises the classical Hamiltonian formulations \cite{Basdevant2007} to the case of Partial Differential Equations (PDEs) such that the space and time variables are treated on the equal footing \cite{Bridges1997} (see also \cite[Chapter 12]{Leimkuhler2004}).

\begin{figure}
  \centering
  \includegraphics[width=0.65\textwidth]{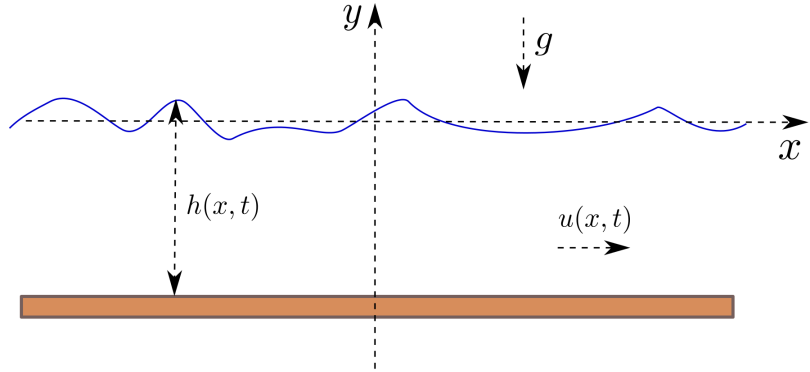}
  \caption{\small\em Sketch of the fluid domain.}
  \label{fig:sketch}
\end{figure}

Let us recall some basic facts about the Serre equations in 2D (one horizontal dimension). Assuming that derivatives are `small' (\ie, long waves in shallow water) but finite amplitudes,\footnote{It should be noted that the steady version of these equations were derived earlier by Rayleigh 
\cite{LordRayleigh1876}.} Serre \cite{Serre1953a, Serre1953} derived the system of equations
\begin{align}
  h_t\ +\ \partial_x\!\left[\,h\,u\,\right]\, &=\ 0,  \label{eqmasse} \\
  u_t\ +\ u\,u_x\ +\ g\,h_x\ +\ \third\,h^{-1}\,\partial_x\!\left[\, h^2\,\gamma\,\right]\,&=\ 0, \label{eqqdmnoncons}
\end{align}
where 
\begin{equation}\label{defgamma}
  \gamma\ =\ h\left[\,u_x^{\,2}\,-\,u_{xt}\,-\,u\,u_{xx}\,\right],
\end{equation}
is the fluid vertical acceleration at the free surface. In these equations, $x$ is the horizontal coordinate, $t$ is the time, $u$ is the depth-averaged horizontal velocity and $h$ is the total water depth (bottom to free surface). A sketch of the domain is shown in Figure~\ref{fig:sketch}.

These equations were independently rediscovered later by Su and Gardner \cite{Su1969} and again by Green, Laws and Naghdi \cite{Green1974}. These approximations being valid in shallow water without assuming small amplitude waves, they are sometimes called \emph{weakly-dispersive fully-nonlinear approximation} \cite{Wu2001a} and are a generalisation of the Saint-Venant \cite{Stoker1957,Wehausen1960} and of the Boussinesq equations \cite{BCS, DMII}. The derivations above are straightforward, one can refer to \cite{Lannes2009}, for example.

From the equations \eqref{eqmasse} to \eqref{defgamma}, one can derive an equation for the conservation of the momentum flux
\begin{align}\label{eqqdmflux}
\partial_t\!\left[\,h\,u\,\right]\, +\ \partial_x\!\left[\,h\,u^2\,+\,\half\,g\,h^2\,+\,\third\,h^2\,\gamma\,
\right]\, &=\ 0, 
\end{align}
and an equation for the conservation of the energy
\begin{align}\label{eqenese}
\partial_t\!\left[\,\half\,h\,u^2\,+\,\sixth\,h^3\/u_x^{\,2}\,+\,\half\,g\,h^2\,\right]\, +\ \partial_x\!
\left[\,(\half\,u^2\,+\,\sixth\,h^2\,u_x^{\,2}+\,g\,h\,+\,\third\,h\,\gamma\,)\,h\,u\,\right]\, &=\ 0,
\end{align} 
as well as an equation for the conservation of the tangential momentum at the free surface
\begin{align}\label{eqqdmtan}
\partial_t\!\left[\,u\,-\,\third\,h^{-1}(h^3\,u_x)_x\,\right]\, +\, \partial_x\!\left[\,\half\,u^2\,+\,g\,h\,
-\,\half\,h^2\,u_x^{\,2}\,-\,\third\,u\,h^{-1}(h^3\,u_x)_x\,\right]\, =\ 0.
\end{align}
Below, we discuss the connection between these conservation laws with the underlying multi-symplectic structure.

The present manuscript is organised as follows. The multi-symplectic formulation is described in the section \ref{sec:ms} and the underlying conservation laws are discussed in the section \ref{sec:cl}. Finally, the main conclusions of this study and some perspective are outlined in the section \ref{sec:disc}.

%%% ------------------------------------------------------------------------ %%%

\section{Multi-symplectic formulation}
\label{sec:ms}

A system of Partial Differential Equations (PDEs) is said to be multi-symplectic if it can written as a system of first-order equations of the form \cite{Bridges1997,Marsden1998}:
\begin{equation}\label{sympgen}
  \mathbb{M}\scal\vz_t\ +\ \mathbb{K}\scal\vz_x\ =\ \grad_{\!z}\,S(\vz),
\end{equation}
where a dot denotes the contracted (inner) product,  $\vz\in\mathds{R}^n$ is a rank-one tensor (vector) of state variables, $\mathbb{M}\in\mathds{R}^{n\times n}$ and $\mathbb{K}\in\mathds{R}^{n\times n}$ are skew-symmetric constant rank-two tensors (matrices) and $S$ is a smooth rank-zero tensor (scalar) function depending on $\vz$. (We use tensor notations because they give more compact formulae than the matrix formalism \cite{Borisenko1979}.) The function $S$ plays the role of the Hamiltonian functional in classical symplectic formulations \cite{Basdevant2007}. Consequently, $S$ is sometimes called the `\emph{Hamiltonian}' function as well. It should be noted that the matrices $\mathbb{M}$ and $\mathbb{K}$ can be (and often are) degenerated \cite{Bridges2006}.

It turns out that the Serre equations \eqref{eqmasse}--\eqref{eqqdmnoncons} have a multi-symplectic structure with $\vz = h\,\ve_1+\phi\,\ve_2+u\,\ve_3+v\,\ve_4+p\,\ve_5+q\,\ve_6+r\,\ve_7+s\,\ve_8$ ($\ve_i$ unitary basis vectors) and
\begin{align}
\mathbb{M}\ &=\ \ve_1 \otimes \ve_2\ -\ \ve_2 \otimes \ve_1 \ +\ \third\,\ve_1 \otimes \ve_5\ -\  \third\,\ve_5 
\otimes \ve_1, \label{defM}\\
\mathbb{K}\ & =\  \third\,\ve_1 \otimes \ve_7\ -\ \third\,\ve_7 \otimes \ve_1\ -\ \ve_2 \otimes \ve_6\ +\ \ve_6 
\otimes \ve_2, \label{defK}\\
S\ &=\,\left(\/\sixth\,v^2\/-\/\half\,u^2\/-\/\third\,s\,u\,v\/\right)h\ -\ \half\,g\,h^2\ +\ \third\,p
\left(\/u\,s\,-\,v\/\right)\, +\ q \left(\/u\/+\/\third\,s\,v\/\right)\, -\ \third\,r\,s. \label{defS}
\end{align}
Indeed, the substitution of these relations into \eqref{sympgen} yields the equations
\begin{align}
  \phi_t\ +\ \third\,p_t\ +\ \third\,r_x\ &=\ \sixth\,v^2\ -\ \half\,u^2\ -\ \third\,s\,u\,v\ -\ g\,h, \label{Sh}\\
  -h_t\ -\ q_x\ &=\ 0, \label{Sphi}\\
  0\ &=\ q\ -\ h\,u\ +\ \third\/s\left(\,p\,-\,h\,v\,\right), \label{Su}\\
  0\ &=\ \third\left(\,h\,v\,-\,p\,\right)\,+\ \third\,s\left(\,q\,-\,h\,u\,\right),\label{Sv}\\
  -\third\,h_t\ &=\ \third\left(\,s\,u\,-\,v\,\right), \label{Sp}\\
  \phi_x\ &=\ u\ +\ \third\,s\,v,\label{Sq}\\
  -\third\,h_x\ &=\ -\third\,s, \label{Sr}\\
  0\ &=\ \third\left(\,p\,u\,+\,q\,v\,-\,r\,-\,h\,u\,v\,\right). \label{Ss}
\end{align}
These equations have the following physical meaning. Equation \eqref{Sr} gives $s = h_x$ so $s$ is the surface slope. Equations \eqref{Su} and \eqref{Sv} yield $p = h v$ and $q = h u$ that are the vertical and horizontal momenta, respectively. It follows that \eqref{Sphi} is the mass conservation $h_t + [h u]_x = 0$ and \eqref{Sp} is the impermeability of the free surface $h_t + u h_x = v$ ($v$ is then the vertical velocity at the free surface). Equation \eqref{Sq} shows that the velocity field is not exactly irrotational for the Serre equations (a well-known result). The definition above of $p$ and $q$ substituted  into \eqref{Ss} gives $r = h u v$. Finally, substituting all the preceding results into \eqref{Sh}, after some algebra, one obtains
\begin{align*}
  \phi_t\ +\ \half\,u^2\ +\ \sixth\,h^2\,u_x^{\,2}\ +\ g\,h\ -\ \third\,h^2\,u_{xt}\ -\ \third\,h\,u\,\partial_x\!\left[\,h\,u_x\,\right]\, =\ 0.
\end{align*}
Differentiating this equation with respect of $x$, eliminating $\phi$ using \eqref{Sq} and exploiting the mass conservation, one gets the equation \eqref{eqqdmtan}.

It should be noted that eliminating $p$, $q$ and $r$, the `\emph{Hamiltonian}' $S$ becomes
\begin{equation}
S\ =\ \half\,h\,u^2\ -\ \sixth\,h\,v^2\ -\ \half\,g\,h^2\ =\ \half\,h\,u^2\ -\ \sixth\,h^3\,u_x^{\,2}\ 
-\ \half\,g\,h^2,
\end{equation}
so $S$ is not a density of total energy, nor a Lagrangian density.

%%% ------------------------------------------------------------------------ %%%

\section{Conservation laws}
\label{sec:cl}

%The main property of the multi-symplecticity is that the multi-symplectic form $\omega$ is conserved, \ie
%\begin{equation*}
%  \partial_t\,\omega\ +\ \partial_x\,\kappa\ =\ 0,
%\end{equation*}
%where $\omega$ and $\kappa$ are the differential two-forms\footnote{$\ud(\vx)$ denotes the differential one-form acting on the vector $\vx$.} 
%\begin{equation*}
%  \omega\ =\ \half\,\ud(\vz)\wedge\ud\!\left(\mathbb{M}\scal\vz\right), \qquad
%  \kappa\ =\ \half\,\ud(\vz)\wedge\ud\!\left(\mathbb{K}\scal\vz\right),
%\end{equation*}
%where $\wedge$ is the exterior product.

A multi-symplectic system of partial differential equations has local conservation laws for the energy and momentum
\begin{equation*}
  \partial_tE(\vz)\ +\ \partial_xF(\vz)\ =\ 0, \qquad
  \partial_tI(\vz)\ +\ \partial_xG(\vz)\ =\ 0, 
\end{equation*}
where 
\begin{gather*}
  E(\vz)\ =\ S(\vz)\ +\ \half\,\vz_x\scal\mathbb{K}\scal\vz, \qquad
  F(\vz)\ =\ -\,\half\,\vz_t\scal\mathbb{K}\scal\vz, \\
  G(\vz)\ =\ S(\vz)\ +\ \half\,\vz_t\scal\mathbb{M}\scal\vz, \qquad
  I(\vz)\ =\ -\,\half\,\vz_x\scal\mathbb{M}\scal\vz.
\end{gather*}
For the Serre equations, from the results of the previous section, we have
\begin{align*}
E\ &=\ \sixth\,r\,h_x\ -\ \sixth\,h\,r_x\ +\ \half\,\phi\,q_x\ -\ \half\,q\,\phi_x\ -\ \half\,g\,h^2\ 
+\ \half\,h\,u^2\ -\ \sixth\,h\,v^2, \\
F\ &=\ \sixth\,h\,r_t\ -\ \sixth\,r\,h_t\ +\ \half\,q\,\phi_t\ -\ \half\,\phi\,q_t, \\
G\ &=\ \sixth\,p\,h_t\ -\ \sixth\,h\,p_t\ +\ \half\,\phi\,h_t\ -\ \half\,h\,\phi_t\ -\ \half\,g\,h^2\ 
+\ \half\,h\,u^2\ -\ \sixth\,h\,v^2, \\
I\ &=\ \sixth\,h\,p_x\ -\ \sixth\,p\,h_x\ +\ \half\,h\,\phi_x\ -\ \half\,\phi\,h_x,
\end{align*}
and using the relations \eqref{Sh}--\eqref{Ss}, after some algebra, one gets the expression of quantities $E$, $F$, $G$ and $I$ in initial physical variables
\begin{align*}
-\/E\ &=\ \half\,h\,u^2\ +\ \half\,g\,h^2\ +\ \sixth\,h^2\,u_x^{\,2}\ -\ 
\partial_x\!\left[\,\half\,\phi\,h\,u\,+\,\sixth\,h^3\,u\,u_x\,\right],\\
-\/F\ &=\,\left(\,\half\,u^2\,+\,\sixth\,h^2\,u_x^{\,2}+\,g\,h\,+\,\third\,h\,\gamma\,\right)h\,u\ +\ 
\partial_t\!\left[\,\half\,\phi\,h\,u\,+\,\sixth\,h^3\,u\,u_x\,\right], \\
G\ &=\ h\,u^2\ +\ \half\,g\,h^2\ +\ \third\,h^2\,\gamma\ +\ 
\partial_t\!\left[\,\half\,\phi\,h\,+\,\sixth\,h^3\,u_x\,\right], \\
I\ &=\ h\,u\ -\ \partial_x\!\left[\,\half\,\phi\,h\,+\,\sixth\,h^3\,u_x\,\right].
\end{align*}
So the momentum and energy conservation equations \eqref{eqqdmflux} and \eqref{eqenese} are recovered, though $-E$, $-F$, $G$ and $I$ are not exactly the densities of energy, energy flux, momentum flux and impulse, respectively.

%%% ------------------------------------------------------------------------ %%%

\section{Discussion}
\label{sec:disc}

In the present manuscript, we discussed the multi-symplectic structure for the Serre equations, which is a very popular model nowadays for long waves in shallow waters. To our knowledge it is the first time that such a structure is reported in the literature. A non-canonical Hamiltonian structure of the Serre equations can be found, for example, in \cite{Johnson2002}. However, we find that the corresponding multi-symplectic structure is simpler and more natural for these equations. Moreover, it allows to treat on the equal footing the space and time variables \cite{Marsden1998}. The advantages of this formulation are well-known \cite{Bridges1997}. 

The multi-symplectic structure of the exact water wave equations being already known \cite{Bridges1997}, it seems natural that approximate equations also have a  multi-symplectic structure. However, Serre's equations being not exactly irrotational, it is not obvious {\em a priori\/} that such a multi-symplectic structure should indeed exists. It is not at all trivial to obtain this structure directly from the Serre equations \eqref{eqmasse}--\eqref{defgamma}. In order to derive the multi-symplectic formulation of the Serre equations, we started form the relaxed variational principle (generalised Hamilton principle) \cite{Clamond2009}. The derivation is then quite transparent, as shown in the appendix. 
   
Serre's equations can be extended in 3D in several ways. One extension of special interest is the so-called {\em irrotational Green--Naghdi equations} \cite{Clamond2009, Kim2001} for which a multi-symplectic structure can be easily obtained following the same route as for the Serre equations, \ie, starting from the relaxed variational principle.

Our study opens some new perspectives to construct structure-preserving integrators for the Serre equations. To our knowledge this research direction is essentially open nowadays. There are some attempts to solve these equations with conventional finite volume \cite{ChazelLannes2010}, pseudo-spectral \cite{Dutykh2011a} and finite element \cite{Mitsotakis2014} methods. However, all these attempts do not guarantee the preservation of the variational (symplectic or multi-symplectic) structures at the discrete level as well. Using the findings reported in this manuscript, one should be able to construct relatively easily finite difference \cite{Ascher2005, Bridges2001, Moore2003a, Wang2003} and pseudo-spectral \cite{Chen2011, Gong2014} schemes, which preserve \emph{exactly} the multi-symplectic conservation law on the discrete level. A numerical comparison of symplectic, multi-symplectic and pseudo-spectral schemes was performed in \cite{Dutykh2013a} on the example of the celebrated 
Korteweg--deVries equation.

%%% ------------------------------------------------------------------------ %%%

\appendix
\section{The workflow pattern}

Our study would not be complete if we did not explain how we arrived to the multi-symplectic structure \eqref{sympgen} of the Serre equations. It is not so trivial to see how this structure appears from equations \eqref{eqmasse}, \eqref{eqqdmnoncons}. However, when we derive the Serre system from the relaxed variational principle \cite{Clamond2009}, a more suitable form of the equations appears. Namely, the relaxed Lagrangian \cite{Clamond2009} under the shallow water ansatz reads (see also \cite{Dutykh2011a})
\begin{equation}
  \mathscr{L}\ =\ (h_t + \vmum\, h_x)\,\phim\ -\ \half\,g\,h^2\ +\ h\left[\,\vmum\vum\ -\ {\half}\,\vum^2\ +\ \third\,\nus\,\vs\ -\ \sixth\,\vs^2\, +\ \phim\,\vmum_x\,\right],
\end{equation}
where $\vmum$, $\nus$ are the Lagrange multipliers. An additional constraint of the free surface impermeability is imposed:
\begin{equation}
  \nus\ =\ h_t\ +\ \vmum\, h_x.
\end{equation}
The corresponding Euler--Lagrange equations are
\begin{eqnarray}
  \udelta\,\vum\/:&& 0\ =\ \vmum\ \ -\ \vum, \label{minL3sw2e1} \\
  \udelta\,\vs\/:&& 0\ =\ h_t\ +\ \vmum h_x\ - \ \vs, \label{minL3sw2e2} \\
  \udelta\,\vmum\/:&& 0\ =\ \vum\ +\ \third\,\vs\, h_x\ -\ \phim_x, \label{minL3sw2e3} \\
  \udelta\,\phim\/:&& 0\ =\ h_t\ +\ \left[\,h\,\vmum\,\right]_x, \label{minL3sw2e4} \\
  \udelta\,h\/:&& 0\ =\ \vmum\vum\ -\ {\half}\,\vum^2\ -\ \sixth\,\vs^2\ -\ \vmum\phim_x\ -\ \phim_t\ -\ g\,h\nonumber\\
  &&\qquad\ -\ \third\,h\left[\,\vs_t\,+\,\vmum\vs_x\,+\,\vs\,\vmum_x\,\right].\label{minL3sw2e5}
\end{eqnarray}
After eliminating $\vmum$ from equations \eqref{minL3sw2e2}--\eqref{minL3sw2e5} thanks to \eqref{minL3sw2e1} and introducing the extra variables $p = h v$, $q = h u$, $r = h u v$ and $s = h_x$, one almost obtains the required system \eqref{Sh}--\eqref{Ss} for the multi-symplectic formulation.

%%% ------------------------------------------------------------------------ %%%

\subsection*{Acknowledgments}
\addcontentsline{toc}{subsection}{Acknowledgments}

D.~\textsc{Clamond} \& D.~\textsc{Dutykh} would like to acknowledge the support of CNRS under the PEPS InPhyNiTi 2015 project FARA.
\smallskip

%%% ------------------------------------------------------------------------ %%%

%%% Bibliography
\addcontentsline{toc}{section}{References}
\bibliographystyle{abbrv}
\bibliography{biblio}

\end{document}